\documentclass[12pt]{article}

\usepackage[T1]{fontenc}
\usepackage{lmodern}
\usepackage{textcomp}

\usepackage{amsmath}%
\usepackage{amssymb}%
\usepackage[table]{xcolor}%

\usepackage{todonotes}
\usepackage[in]{fullpage}%

\usepackage[amsmath,thmmarks]{ntheorem}%

\theoremseparator{.}%

\usepackage{titlesec}%
\titlelabel{\thetitle. }%
\usepackage{xcolor}%
\usepackage{mleftright}%
\usepackage{xspace}%
\usepackage{graphicx}
\usepackage{hyperref}%

\usepackage{hyperref}%
\hypersetup{%
      unicode,
      breaklinks,%
      colorlinks=true,%
      urlcolor=[rgb]{0.25,0.0,0.0},%
      linkcolor=[rgb]{0.5,0.0,0.0},%
      citecolor=[rgb]{0,0.2,0.445},%
      filecolor=[rgb]{0,0,0.4},
      anchorcolor=[rgb]={0.0,0.1,0.2}%
}
\usepackage[ocgcolorlinks]{ocgx2}

\usepackage{lmodern}

\theoremseparator{.}%

\theoremstyle{plain}%
\newtheorem{theorem}{Theorem}[section]

\newtheorem{lemma}[theorem]{Lemma}

\newtheorem{fact}[theorem]{Fact}

\newtheorem{prop}[theorem]{Proposition}

\newtheorem{remark}[theorem]{Remark}

\newtheorem*{remark:unnumbered}[theorem]{Remark}%

\newcommand{\myqedsymbol}{\rule{2mm}{2mm}}

\theoremheaderfont{\em}%
\theorembodyfont{\upshape}%
\theoremstyle{nonumberplain}%
\theoremseparator{}%
\theoremsymbol{\myqedsymbol}%
\newtheorem{proof}{Proof:}

\newenvironment{proofof}[1]
{\par\noindent\textit{Proof of #1:}}
{\hfill$\myqedsymbol$\par}

\providecommand{\emphind}[1]{}%
\renewcommand{\emphind}[1]{\emph{#1}\index{#1}}

\definecolor{blue25emph}{rgb}{0, 0, 11}

\providecommand{\emphic}[2]{}
\renewcommand{\emphic}[2]{\textcolor{blue25emph}{%
      \textbf{\emph{#1}}}\index{#2}}

\providecommand{\emphi}[1]{}%
\renewcommand{\emphi}[1]{\emphic{#1}{#1}}

\definecolor{almostblack}{rgb}{0, 0, 0.3}

\providecommand{\emphw}[1]{}%
\renewcommand{\emphw}[1]{{\textcolor{almostblack}{\emph{#1}}}}%

\providecommand{\emphOnly}[1]{}%
\renewcommand{\emphOnly}[1]{\emph{\textcolor{blue25}{\textbf{#1}}}}

\newcommand{\HLink}[2]{\hyperref[#2]{#1~\ref*{#2}}}
\newcommand{\HLinkSuffix}[3]{\hyperref[#2]{#1\ref*{#2}{#3}}}

\providecommand{\deflab}[1]{}
\renewcommand{\deflab}[1]{\label{def:#1}}

\providecommand{\eqlab}[1]{}%
\renewcommand{\eqlab}[1]{\label{equation:#1}}

\newcommand{\remove}[1]{}%

\newcommand{\R}{\mathbb{R}}

\usepackage[inline]{enumitem}

\newlist{compactenumA}{enumerate}{5}%
\setlist[compactenumA]{topsep=0pt,itemsep=-1ex,partopsep=1ex,parsep=1ex,%
   label=(\Alph*)}%

\newlist{compactenuma}{enumerate}{5}%
\setlist[compactenuma]{topsep=0pt,itemsep=-1ex,partopsep=1ex,parsep=1ex,%
   label=(\alph*)}%

\newlist{compactenumI}{enumerate}{5}%
\setlist[compactenumI]{topsep=0pt,itemsep=-1ex,partopsep=1ex,parsep=1ex,%
   label=(\Roman*)}%

\newlist{compactenumi}{enumerate}{5}%
\setlist[compactenumi]{topsep=0pt,itemsep=-1ex,partopsep=1ex,parsep=1ex,%
   label=(\roman*)}%

\newlist{compactitem}{itemize}{5}%
\setlist[compactitem]{topsep=0pt,itemsep=-1ex,partopsep=1ex,parsep=1ex,%
   label=\ensuremath{\bullet}}%

\numberwithin{equation}{section}%

\begin{document}

\title{Polyak-Łojasiewicz inequality is essentially no more general than strong convexity for $C^2$ functions}

\author{%
    Aziz Ben Nejma%
   \thanks{Université Paris Cité and Sorbonne Université, CNRS, Laboratoire de Probabilités, Statistique et Modélisation, F-75013 Paris, France and DMA, Ecole normale supérieure, Université
PSL, CNRS, 75005 Paris, France.}%
}

\maketitle
\begin{abstract}
    The Polyak-Łojasiewicz (PŁ) inequality extends the favorable optimization properties of strongly convex functions to a broader class of functions. In this paper, we prove a theorem (also obtained by Criscitiello, Rebjock and Boumal in an earlier blog post) showing that the richness of the class of PŁ functions is rooted in the nonsmooth case since sufficient regularity forces them to be essentially strongly convex. More precisely, if $f$ is a $C^2$ PŁ function having a bounded set of minimizers, then it has a unique minimizer and is strongly convex on a sublevel set of the form $\{f\leq a\}$. We show that this implies a result of Asplund on properties of the squared distance function, and discuss some consequences on smoothness assumptions in results in the literature.
\end{abstract}

\setcounter{secnumdepth}{3}

\section{Introduction}

Let $f$ denote a $C^2$ function defined on $\mathbb{R}^n$, the function $f$ is said to satisfy a Polyak-Łojasiewicz inequality (or equivalently to be PŁ) if there exists a positive constant $C$ such that 
$$\lVert \nabla f \rVert^2 \geq C (f-\mathrm{inf}f).$$
The Polyak-Łojasiewicz inequality is implied by strong convexity, and although being a significantly weaker condition, it is as handy when establishing the convergence of gradient descent. More precisely, the function $f$ satisfies a Polyak-Łojasiewicz inequality with constant $C$ if and only if $$f(y_t)-\inf f \leq e^{-Ct} (f(y_0)-\inf f)$$ where $(y_t)_{t\geq 0}$ is the gradient flow starting from $y_0 \in \R^d$, defined by $\dot{y_t}=-\nabla f(y_t)$.
This fact was first observed by Polyak \cite{polyak1963gradient}. And while being weaker than most assumptions that guarantee linear convergence for several optimization methods, Polyak-Łojasiewicz inequalities still yield the same convergence bounds for continuous as well as discrete optimization methods in the nonconvex setting \cite{polyak1963gradient,Attouch2009, Attouch2013, karimi2016linear}.
These inequalities, or local versions of them, are easily satisfied and cover a wider range of functions when compared to strong convexity :
they do not only hold in convex settings, in fact, Łojasiewicz then Kurdyka \cite{lojasiewicz1963propriete, kurdyka1998gradients} showed Łojasiewicz inequalities in o-minimal structures which include analytic and subanalytic functions and these functions therefore satisfy PŁ-type inequalities around their minimizers, i.e.
$$\lVert \nabla f(x) \rVert^q \geq C (f(x)-\underset{\mathbb{R}^n}{\mathrm{min}}f)$$
for some $q>1$ and for $x$ close enough to $\mathrm{argmin}(f)$.

Moreover, Polyak-Łojasiewicz inequalities and most of their properties extend to nonsmooth settings through limiting subdifferentials or using the framework of gradient flows in metric spaces. As a result, the Polyak-Łojasiewicz inequality became a standard assumption in many problems (see, for example, \cite{vaswani2019fast} or \cite{chewi2024ballistic}) and an object of study in its own right.

A significant and prototypical class of PŁ functions of interest to us consists of functions of the form
$$x \longmapsto \mathrm{dist}^2(x,F)$$
where $F \subseteq \R^d$ is any closed set and $\mathrm{dist}$ is the usual distance function. A more general version of this statement can be found, for instance, in \cite{garrigos2023square}.

\noindent Finally, it is worth emphasizing that unlike strict convexity or $\gamma$-strong convexity for a given $\gamma >0$, which have a local description, the PŁ condition with constant $C>0$ is defined as a global property, which may involve additional subtleties.

\section{Main result}
The main result of the present paper describes the possible sets of minimizers of a $C^2$ Polyak-Łojasiewicz function. In many respects, $C^{1,1}_{loc}$ and $C^2$ functions behave almost identically. For instance, $C^{1,1}$ functions are twice differentiable almost everywhere. However, $C^2$ Polyak-Łojasiewicz functions are characterized by a surprising rigidity in their optimization landscape. Here is the main theorem.
\begin{theorem}\label{thm1}
    Let $f \in C^2(\mathbb{R}^n)$ and assume that $f$ satisfies a Polyak-Łojasiewicz inequality and that $\mathrm{argmin}(f)$ is a bounded set. Then $\mathrm{argmin}(f)$ consists of a single element. 
\end{theorem}
\begin{remark}
\normalfont
    After the initial public posting of this paper, we became aware that the result of Theorem \ref{thm1} had previously been established in a blog post by Criscitiello, Rebjock and Boumal \cite{criscitiello2025} based on earlier work \cite{Rebjock2025}, with a proof similar to ours.
\end{remark}

\begin{remark}
\normalfont
Strongly convex functions always having a unique minimizer, Theorem \ref{thm1} can  therefore be restated as follows :

\noindent \textit{Let $f \in C^2(\mathbb{R}^n)$ and assume that $f$ satisfies a Polyak-Łojasiewicz inequality and that $\mathrm{argmin}(f)$ is a bounded set. Then there exists an open convex set $\Omega$ containing $\mathrm{argmin}(f)$ such that $f_{|\Omega}$ is a strongly convex function.}
\end{remark}

\begin{remark}
\normalfont
    The condition $f \in C^2(\mathbb{R}^n)$ is essential in Theorem \ref{thm1}.
    Take \[f : x \in \mathbb{R} \longmapsto \frac{1}{2}\mathrm{dist}(x,[0,1])^2. \] In this case, it holds that $f \in C^{1,1}_{loc}(\mathbb{R}^n)$ and that it satisfies a global PŁ inequality yet the set of its critical points is compact without being a singleton. 
\end{remark}

\begin{remark}\label{rqp}
\normalfont
    The compactness of $\mathrm{argmin}(f)$ is also essential. In fact, any positive semidefinite quadratic form is a smooth PŁ function, yet the set of its minimizers is a linear subspace of $\mathbb{R}^n$.
    More generally, let $f$ be a PŁ function defined on $\mathbb{R}^m$. We can define the function $$\Tilde{f} : (x,x') \in  \mathbb{R}^m \times \mathbb{R}^k \longmapsto f(x).$$
    The new function $\Tilde{f}$ has same regularity as $f$ and satisfies a PŁ inequality with same constant, however, $\mathrm{argmin}(\tilde{f})= \mathrm{argmin}(f) \times \mathbb{R}^k$. 
\end{remark}
\begin{remark}
\normalfont Theorem \ref{thm1} does not assume that $\mathrm{argmin}(f) \neq \emptyset$. In fact, for all PŁ functions $f \in C^{1,1}_{loc}(\R^d)$, it holds that $\mathrm{argmin}f \neq \emptyset$. This can be seen as a consequence of Theorem \ref{length}, as gradient flow trajectories of PŁ functions have finite length.
\end{remark}
\begin{remark}
\normalfont
As already recalled, one important class of PŁ functions consists of all functions of the form $x \longmapsto \mathrm{dist}(x,F)^2$, where $F$ is a (closed) subset of $\mathbb{R}^n$. Theorem \ref{C2dist} implies that despite all these functions being PŁ, the only ones being $C^2$ are those such that the set of their minimizers $F$ is an affine subspace of $\mathbb{R}^n$.
\noindent Having that in mind (along with Theorem \ref{thm1} and remark \ref{rqp}), it may be tempting to think that the set of minimizers of a $C^2$ PŁ function is an affine subspace of $\mathbb{R}^n$, but as one would reasonably expect, this is not true !

\noindent Take any function $g \in C^2(\mathbb{R}^n)$ and define 
$$f : (x,y) \in \mathbb{R}^n \times \mathbb{R} \longmapsto (y-g(x))^2.$$
The function $f$ is $C^2$ (it can be made arbitrarily smooth) and satisfies a PŁ inequality. However, in general, $\mathrm{argmin}(f)$ is not an affine subspace.
\end{remark}

\begin{remark}
\normalfont
    There is no single general reason why PŁ functions almost always fail to be $C^2$. For instance, in the case of $\mathrm{dist}(.,[0,1])^2$, the issue lies on the boundary of the minimizing set, whereas in the case of $\mathrm{dist(.,\mathbb{S}^1)}^2$, the issue lies far from the minimizing set, being the non-differentiability at $0$.
\end{remark}

\begin{remark}
\normalfont
    The use of the word \textit{essentially} in the title is intended in the sense of topology, not in that of optimization. In particular, Theorem \ref{thm1} does not relate global optimization properties of $C^2$ PŁ functions to those of strongly convex functions as it only provides local information around minimizers. Instead, it shows that these functions are similar in the sense that their level sets share the same topological properties.
\end{remark}
 
\noindent As a corollary, we recover a classical result of Asplund \cite{Asplund1969, PoliquinRockafellarThibault2000,doi:10.1142/12797} on the connection between the regularity of a closed subset of $\R^n$ and the regularity of the distance function to it. More precisely,
\begin{theorem}\label{C2dist}
    Let $F$ be any closed subset of $\mathbb{R}^n$. Then 
    \begin{enumerate}[]
        \item $d_F : x \mapsto \mathrm{dist}^2(x,F)$ is $C^{1,1}_{loc}$ if and only if $F$ is convex.
        \item $d_F : x \mapsto \mathrm{dist}^2(x,F)$ is $C^2$ if and only if $F$ is an affine subspace of $\mathbb{R}^n$.
    \end{enumerate}
\end{theorem}

\medskip
To our best knowledge, although several optimization methods have been thoroughly investigated under a Polyak-Łojasiewicz condition, the topological properties of level sets of PŁ functions and in particular properties of their minimal sets were not studied. Yet, several theorems require assumptions on PŁ or locally PŁ functions' minimizers topology like connectedness, or being a discrete set. This has often led to unnecessary assumptions such as uniqueness of a minimizer of a smooth proper PŁ function in \cite{chewi2024ballistic}.

\paragraph{Acknowledgements.} I am deeply grateful to Max Fathi, under whose supervision I carried out my internship, for his guidance, support and very helpful suggestions. I also thank Guillaume Garrigos and Austin Stromme for fruitful discussions.\\
The author has received support under the program "Investissement d'Avenir" launched by the French Government and implemented by ANR, with the reference ANR-18-IdEx-0001 as part of its program "Emergence".  He was also supported by the Agence Nationale de la Recherche (ANR) Grant ANR-23-CE40-0003 (Project CONVIVIALITY).

\section{Outline of Proof of Theorem \ref{thm1}}

In the sequel, $f$ denotes a $C^2$ function defined on $\mathbb{R}^n$ and satisfying a Polyak-Łojasiewicz inequality with constant $C$, i.e.
$$\lVert \nabla f \rVert^2 \geq C (f-\inf(f)).$$
Without loss of generality, we may assume $C=1$ and $\inf(f)=0$.\\
Let $K$ denote the set of critical points of $f$, since $f$ is PŁ, $K=\mathrm{argmin(f)}$ so that it is supposed compact, and is nonempty by Remark \ref{fin}.

\begin{proof}

The proof of Theorem \ref{thm1} relies on the following lemmas.

\begin{lemma}\label{l2}
$K$ is a contractible set, that is, $K$ has the homotopy type of a point.
\end{lemma}

\begin{lemma}\label{l3}
    K is a closed submanifold of $\mathbb{R}^n$, that is, $K$ is a compact submanifold of $\mathbb{R}^n$ with no boundary.
\end{lemma}

\begin{remark}\label{rk3}
\normalfont
    The results of lemmas \ref{l2} and \ref{l3} (except compactness in Lemma \ref{l3}) remain true without any assumptions on $K=\mathrm{argmin}(f)$ whenever $f \in C^2(\mathbb{R}^n)$ and the restriction of $f$ to any compact subset of $\mathbb{R}^n$ satisfies a PŁ inequality, i.e. $$\forall R>0~\exists C_R>0~\forall x \in \mathbb{R}^n~\big( \lVert x \rVert < R \implies \lVert \nabla f(x) \rVert^2 \geq C_R(f(x)-\underset{\mathbb{R}^n}{\inf}f) \big).$$
\end{remark}
Let us admit for now both of the lemmas above. Then $K$ is a contractible closed manifold. Let $d$ be its dimension, then $d=0$. In fact, for a closed connected $d$-dimensional manifold $M$ having dimension $d\geq1$, its $d^{{th}}$ homology group\footnote{The only use of homology and algebraic topology arguments in the paper occurs here. For the reader's convenience we recall the necessary facts : given a manifold, one can associate a sequence of abelian groups, the homology groups, denoted $H_0,H_1,\dots$.
    The homology of a topological space is invariant under homotopy, that is, under continuous deformation of the space.\\ Moreover, contractible spaces are topological spaces that have the homotopy type of a point. A continuous image of a contractible space is contractible. For a contractible space, as for a point, $H_n \cong \{0\}$ for $n\geq 1$.} is nontrivial as we have (see for instance \cite[Theorem 3.26]{MR1867354}) $$H_d(M,\mathbb{Z}/2\mathbb{Z}) \cong \mathbb{Z}/2\mathbb{Z} $$
which cannot hold for a contractible set.
Thus, $K$ is a $0$-dimensional connected manifold, that is, a single point.
\end{proof}

\section{Proof of the lemmas}

In order to prove the lemmas, we recall some useful well-known facts about PŁ functions. For proofs, the reader is referred to \cite{bolte2010characterizations} and \cite{karimi2016linear}.

\subsection{A few basic facts about PŁ functions}

\begin{theorem}\label{cv}
Suppose $f \in C^2(\mathbb{R}^n)$ satisfies $\lVert \nabla f \rVert^2 \geq C (f-\inf(f))$. Let $y$ be the solution of the ODE
\[
\left \{
\begin{array}{c @{=} c}
    \dot{y}(t) & -\nabla f (y(t)) \\
    y(0) & y_0
\end{array}
\right.
\]
Then $f(y(t))-\inf(f) \leq e^{-Ct}(f(y_0)-\inf(f))$.
\end{theorem}

\begin{theorem}\label{length}
Suppose $f \in C^{1,1}_{loc}(\mathbb{R}^n)$ satisfies $\lVert \nabla f \rVert^2 \geq f$ and that $\inf(f)=0$. Let $y$ be the solution of the ODE
\[
\left \{
\begin{array}{c @{=} c}
    \dot{y}(t) & -\nabla f (y(t)) \\
    y(0) & y_0
\end{array}
\right.
\]
Then $$\int_0^{+\infty} \lVert  \dot{y}(t)\rVert dt \leq  2\sqrt{f(y_0)}.$$

\end{theorem}

\begin{remark}\label{fin}
\normalfont Theorem \ref{length} implies that gradient flow trajectories have finite lengths and therefore $\mathrm{argmin(f)\neq \emptyset}$ \cite{bolte2010characterizations}.
\end{remark}

\begin{theorem}\label{growth}
    Suppose $f \in C^2(\mathbb{R}^n)$ satisfies $\lVert \nabla f \rVert^2 \geq f$ and that $\mathrm{min}(f)=0$. \\Then $f(x) \geq \mathrm{dist}(x,\mathrm{argmin}(f))^2$ for all $x \in \mathbb{R}^n$.
\end{theorem}

\noindent In what follows, we keep the same assumptions and notations as in the previous section.

\subsection{Proof of Lemma \ref{l2}}

\begin{prop}\label{propo21}
    For $t\geq 0$ and $x \in \mathbb{R}^n$, let $\varphi(t,x)$ be the flow of $(-\nabla f)$, that is, $\varphi(t,x)$ is the value at time $t$ of the the solution of $\dot{y} = -\nabla f(y)$ with initial condition $y(0)=x$.
    Then $\varphi(t,.)$ converges uniformly on compact sets as $t \to +\infty$.
\end{prop}

\begin{proof}
Our proof closely follows \cite{bolte2010characterizations}. Since $f$ is $C^2$, $\varphi$ is well-defined by Cauchy-Lipschitz 
    and for every $t \geq 0$, $\varphi(t,.)$ is a $C^1$ diffeomorphism of $\mathbb{R}^n$ to itself.
Moreover, Theorem \ref{length} states that the curve $t\mapsto \varphi(t,x)$ has finite length. More precisely, for every $x \in \mathbb{R}^n$ the following inequality holds.
$$ \int_0^{+\infty} \lVert  \frac{d}{dt} \varphi(t,x)\rVert dt \leq  2\sqrt{f(x)}.$$
Hence, for every nondecreasing sequence of nonnegative times $(t_n)_n$ such that $t_n \to +\infty$, it holds that
$$\displaystyle\sum_{k=0}^{+\infty} \lVert \varphi(t_{k+1},x)-\varphi(t_{k},x) \rVert \leq \displaystyle\sum_{k=0}^{+\infty} \lVert\int_{t_k}^{t_{k+1}}   \frac{d}{dt} \varphi(t,x) dt\rVert \leq \int_0^{+\infty} \lVert  \frac{d}{dt} \varphi(t,x)\rVert dt  \leq 2\sqrt{f(x)}.$$
Therefore, $\varphi(t,.)$ converges normally on compact subsets of $\mathbb{R}^n$, which ends the proof.
\end{proof}

Now, thanks to Proposition \ref{propo21}, we may define $\varphi_\infty$ as the locally uniform limit of the functions $\varphi(t,.)$ when $t \to +\infty$.
Since each $\varphi(t,.)$ is continuous, $\varphi_\infty$ is also a continuous mapping from $\mathbb{R}^n$ to $K$. Indeed, $\varphi_\infty$ takes its values in $K$ because of Theorem \ref{cv}.
Moreover, since $\varphi(t,x) = x$ for all $t\geq 0$ and all $x \in K$, the limit function $\varphi_\infty : \mathbb{R}^n \xrightarrow[]{} K$ is a retraction. $K$ is therefore a contractible subset of $\mathbb{R}^n$, which concludes the proof of Lemma \ref{l2}.

\subsection{Proof of Lemma \ref{l3}}
In this subsection, we denote by $H$ the Hessian matrix associated to $f$, that is $H= \nabla^2f$.
Since $K=\mathrm{argmin}(f)$, $H(x)$ is a positive semidefinite matrix for every $x \in K$. We also denote by $B(x,r)$ the open ball of radius $r$ centered at $x$ and by $S(x,r)$ the Euclidean sphere of radius $r$ centered at $x$.

\begin{prop}\label{rank}
    $H$ has constant rank on $K$. 
\end{prop}

\begin{proof}
Since $f$ satisfies a Polyak-Łojasiewicz inequality with constant $C=1$, a Taylor expansion shows that positive eigenvalues of $H(x)$ lay in the interval $[\frac{1}{2}, +\infty)$ for every $x \in K$. Indeed, for $x \in K$, $f(x)=0$ and $\nabla f(x)=0$, therefore, for $u \in \mathbb{R}^n$,
\[
\left \{
\begin{array}{l}
    f(x+tu)=f(x)+\langle \nabla f(x), u\rangle t+\frac{1}{2}\langle H(x)u, u\rangle t^2+ o(t^2)=\frac{1}{2}\langle H(x)u, u \rangle t^2 + o(t^2) \\
    \nabla f(x+tu)=\nabla f(x)+ tH(x)u +o(t)=tH(x)u+o(t)   
\end{array}
\right.
\]

Writing $f(x+tu) \leq \lVert\nabla f(x+tu)\rVert^2$ yields $\frac{1}{2}  \langle H(x)u,u\rangle \leq \rVert H(x)u\rVert^2$, which means that the spectrum of $H$ is  contained in $\{0\} \cup[\frac{1}{2},+\infty)$, as $H$ is positive semidefinite.

\medskip
    Let $m \in \mathbb{N}_{\geq 0}$ denote the minimum rank of $H$ on $K$. Define $E$ as the set of $x \in K$ such that $ \mathrm{rank}(H(x))=m$. $E$ is therefore a closed subset of $K$ by the lower semicontinuity of the rank.

    Moreover, let $(x_k)_{k \in \mathbb{N}}$ be a sequence of elements of $K$ converging to $x \in E$. By continuity of the spectrum, the number of eigenvalues of $H(x_k)$ in the interval $(-\frac{1}{4},\frac{1}{4})$ converges to $n-m$ as $k \to +\infty$. Yet the only eigenvalue $H(x_k)$ may have in $(-\frac{1}{4},\frac{1}{4})$ is $0$. Thus $\mathrm{rank}(H(x_k))=m$ for large values of $k$.
    This shows that $E$ is a nonempty closed-open subset of $K$. Thus $E=K$ since Lemma \ref{l2} implies in particular that $K$ is connected.
\end{proof}

\begin{prop}\label{inj}
    Let $x \in K$. Let $\pi$ be the orthogonal projection on $\mathrm{Ker}(H(x))$.
    Then there exists $r>0$ such that $\pi_{|K \cap B(x,r)}$ is injective.
\end{prop}

\begin{proof}
Let us denote by $A$ the matrix $H(x)$. For $y,z \in \mathbb{R}^n$ such that $\pi(y)=\pi(z)$. 
\begin{equation}
\begin{split}
\langle \nabla f(z)-\nabla f(y),z-y\rangle & = \int_0^1 \langle H\big((1-t)y+tz\big)(z-y),z-y\rangle dt \\
 & = \int_0^1 \langle A(z-y),z-y\rangle + \langle (H_t-A)(z-y),z-y\rangle dt\\
 & \geq \frac{1}{2} \lVert z-y \rVert^2 - \underset{t \in [0,1]}{\sup} \lVert  H_t-A \rVert_{op}\lVert z-y \rVert^2 \\
 & =\big (\frac{1}{2}-\underset{t \in [0,1]}{\sup} \lVert  H_t-A \rVert_{op}\big)\lVert z-y \rVert^2
\end{split}
\nonumber
\end{equation}
where $H_t= H\big((1-t)y+tz\big)$.

The function $f$ being $C^2$, there exists a positive radius $r$ such that $\lVert H(y)-A\rVert_{op} \leq \frac{1}{4}$ for all $y \in B(x,r)$.
Therefore, whenever $y,z \in B(x,r)$ are such that $\pi(y)=\pi(z)$, one cannot have $\nabla f(y) = \nabla f(z)$ ($=0$ when $y,z \in K)$ unless $y=z$.
\end{proof}

\begin{prop}\label{surj}
Let $x \in K$ and $r>0$. Then $\pi \big( K \cap B(x,r)\big)$ contains a neighborhood of $\pi(x)$ in $\mathrm{Ker}(H(x))$, where $\pi$ is the orthogonal projection on $\mathrm{Ker}(H(x))$.
\end{prop}

\begin{proof}
    We may assume, by taking a smaller $r$ if necessary, that $\pi_{|K \cap \overline{B(x,r)}}$ is injective and that $\lVert H(z)-H(x)\rVert_{op} \leq \dfrac{1}{4}$ for all $ z \in B(x,r) \cap K$. Let $v $ be an element of $\mathrm{Ker}(H(z))^{\perp} \cap \mathrm{Ker}(H(x))$ where $z \in K$. It holds that
    $$ \dfrac{\lVert v \rVert}{2} \leq \lVert H(z)v \rVert \leq \lVert H(x)v \rVert +\dfrac{\lVert v \rVert}{4}=\dfrac{\lVert v \rVert}{4}$$
    where, in the first inequality, we used that the spectrum of $H(z)$ is contained in $\{0\} \cup [\frac{1}{2},+\infty)$. Therefore, $v=0$ and $\mathrm{Ker}(H(z))^{\perp} \oplus \mathrm{Ker}(H(x))=\mathbb{R}^n$ since $H(x)$ and $H(z)$ have the same rank by proposition \ref{rank}.

     Since $\pi_{|K \cap \overline{B(x,r)}}$ is injective, $\pi(x) \notin\pi\big(K \cap S(x,r)\big)$ and there exists therefore a positive radius $r'<r$ such that $\pi\big(K \cap S(x,r)\big) \cap B(\pi(x),r') = \emptyset$.

    Assume, for the sake of contradiction, that there is no such neighborhood of $\pi(x)$ as in the statement of the proposition.
    Let $a$ be a point of $B(\pi(x),\frac{r'}{8}) \cap \mathrm{Ker}(H(x))$ such that $a\notin\pi \big( K \cap \overline{B(x,r)}\big)$. Since $\pi \big( K \cap \overline{B(x,r)}\big)$ is closed, there is a smallest $\delta>0$ such that $S(a,\delta) \cap \pi \big( K \cap \overline{B(x,r)}\big) \neq \emptyset$ \big(and $\delta\leq \dfrac{r'}{8}$ \big).

    Let $z \in \overline{B(x,r)} \cap K$ such that $\pi(z) \in S(a,\delta)$. By definition of $r'$, $z\in B(x,r)$. Let $v= \pi(z)-a$. Then there exists $u \in \mathrm{Ker}(H(z))$ such that $\langle u, v\rangle \neq 0$. In fact, since $v \in \mathrm{Ker}(H(x))$,
    $$v^{\perp} \cap  \mathrm{Ker}(H(z)) = \big(\mathrm{Span}(v) \oplus \mathrm{Ker}(H(z))^{\perp}\big)^{\perp}   $$ and this subspace has dimension $$n-1-\mathrm{rank}(H(z)) < n-\mathrm{rank}(H(z)) = \mathrm{dim}(\mathrm{Ker}(H(z))).$$

    Finally, observe that the growth condition (Theorem \ref{growth}) and a Taylor expansion of $f$ at $z$ imply that $$\mathrm{dist}(z+tu,K)^2 \leq f(z+tu) = o(t^2).$$
    Therefore, for every $t \in \mathbb{R}$, there exists a critical point $k_t \in K$ such that $\lVert k_t-z-tu\rVert = o(t)$, that is, $k_t= z+tu + o(t)$.
    Now, write
\begin{equation}
\begin{split}
\lVert \pi(k_t)-a \rVert^2 
&= \lVert \pi(z)-a +t\pi(u) +o(t) \rVert^2 \\
& = \lVert \pi(z)-a \rVert^2 + 2t\langle \pi(u),\pi(z)-a\rangle +o(t) \\
& =\lVert \pi(z)-a \rVert^2 + 2t\langle u,\pi(z)-a\rangle +o(t)\\
& <\delta^2 \text{ for a chosen sign of $t$ with $|t|$ arbitrarily small}.
\end{split}
\nonumber
\end{equation}

The inequality above contradicts the minimality of $\delta$.
Hence, $\pi \big( K \cap B(x,r)\big)$ contains a neighborhood of $\pi(x)$ in $\mathrm{Ker}(H(x))$.
\end{proof}

Combining Proposition \ref{inj} and Proposition \ref{surj} shows that every point of $K$ has a neighborhood that is homeomorphic to the open unit ball of $\mathbb{R}^{n-m}$ where $m$ is the rank of the matrix $H$ at any point of $K$, therefore $K$ is indeed a closed manifold.

\begin{remark}
\normalfont
    It is natural to ask whether Theorem $\ref{thm1}$ remains valid if we assume $f$ to satisfy a PŁ inequality only locally, for example, on a domain containing $\mathrm{argmin(f)}$. The answer is negative. In fact, the square distance to the sphere $f ~\colon~ x\mapsto \mathrm{dist}^2(x,\mathbb{S}^{n-1})$ is PŁ and smooth on $\mathbb{R}^n\backslash\{0\}$, and $f_{|\mathbb{R}^n \backslash B(0,\frac{1}{2})}$ can be extended to a smooth function on $\mathbb{R}^n$ satisfying a PŁ inequality on $\mathbb{R}^n \backslash B(0,\frac{1}{2})$. However, $\mathbb{S}^{n-1}=\mathrm{argmin}(f)$, is compact without being a single point. Note that in this case, only Lemma \ref{l2} fails while Lemma \ref{l3} remains valid under a local PŁ assumption near the minimizers. Finally, the result of Theorem \ref{thm1} remains valid if we assume that $f$ satisfies a PŁ inequality on a contractible domain that is stable under the flow $\varphi(t,.)$, the proof being exactly the same.
\end{remark}

\section{Additional results on the gradient flows of general PŁ functions}

In this section, we keep the same notations but we do not suppose anymore that $\mathrm{argmin}(f)$ is compact, so that it is an embedded submanifold of $\mathbb{R}^n$ by Remark \ref{rk3}. We still assume without loss of generality that $\mathrm{min}(f)=0$ and that $f$ satisfies a PŁ inequality with constant $1$.
We also assume that $\mathrm{dim}(\mathrm{argmin}(f))<n$ since $\mathrm{dim}(\mathrm{argmin}(f))=n$ corresponds to the case where $f$ is constant ($\mathrm{argmin}(f)$ would be a nonempty closed-open subset of $\mathbb{R}^n$).
\begin{theorem}
    Let $x \in \mathrm{argmin}(f)$, then $\{x\} \subsetneq \varphi_\infty^{-1}(\{x\})$.
\end{theorem}

\begin{proof}
We will prove that $\varphi_\infty\big(\mathbb{R}^n\backslash \mathrm{argmin}(f)\big)=\mathrm{argmin}(f)$.

\noindent Let $x \in \mathrm{argmin}(f)$, Let $B$ be a ball centered at $x$ and $S$ its boundary. Let $y \in B\backslash \mathrm{argmin}(f)$, then $y \in \varphi([0,+\infty) \times S)$. In fact, denoting by $\psi$ the flow of $\nabla f$, the curve of $\psi(t,y)$ intersects $S$ for some $t_0 \geq 0$, since $f$ has no critical points except its minimizers. Therefore, we have $y \in \varphi([0,+\infty) \times S)$.

\noindent Moreover, Theorem \ref{length} implies that $\lVert \varphi_\infty(y)-y \rVert \leq 2\sqrt{f(y)}  $. Therefore,
$$\mathrm{dist}(x,\varphi_\infty(S)) \leq \lVert x - \varphi_\infty(y) \rVert \leq \lVert x-y \rVert + 2\sqrt{f(y)}.$$

\noindent Since $\mathrm{dim}(\mathrm{argmin}(f))<n$, we can take the limit as $y \to x$ and deduce that $$\mathrm{dist}(x,\varphi_\infty(S))=0.$$
The result follows from the compactness of $S$ and the continuity of $\varphi_\infty$.
\end{proof}

\section{Regularity of distance functions}
In this section, we recover some well-known results on the regularity of distance functions, taking advantage from the PŁ framework as we prove Theorem \ref{C2dist}, originally due to Asplund \cite{Asplund1969}.

Whenever $F$ is a closed subset of $\mathbb{R}^n$, we define $d_F : x \mapsto \frac{1}{2}\mathrm{dist}^2(x,F)$ to be (half of) the square distance to $F$ function and we denote by $\mathrm{proj}_F$ the projection on $F$, such that for all $x \in \mathbb{R}^n$, $\mathrm{proj}_F(x)$ is the set of points in $F$ that have minimal distance to $x$.
\\
\par If $\mathrm{proj}_F(x)$ is a single point for every $x$, we shall write by abuse of notation $\mathrm{proj}_F(x)=y$ where $y \in  \mathbb{R}^n$ is such that $\{y\}=\mathrm{proj}_F(x)$. We finally recall the following fact (see, e.g., \cite[Lemma 7]{garrigos2023square}):
\begin{fact}
    Let $F$ be a closed subset of $\mathbb{R}^n$, and suppose $d_F$ is $C^1$, then $\mathrm{proj}_F(x)$ is a single point for every $x \in \mathbb{R}^n$ and $\nabla d_F(x) = x-\mathrm{proj}_F(x)$.
\end{fact}

In order to prove Theorem \ref{C2dist}, we fix a closed subset $F \subseteq \mathbb{R}^n$ such that $d_F$ is $C^{1,1}_{loc}$. We denote by $\varphi$ the flow of $-\nabla d_F$ defined as in the previous sections, which is well-defined by the Cauchy-Lipschitz theorem.
The proof relies on the following lemma.

\begin{lemma}\label{Cauchylip}
    $\varphi(t,x)=\mathrm{proj}_F(x)+e^{-t}(x-\mathrm{proj}_F(x))$ and $\mathrm{proj}_F(\varphi(t,x))=\mathrm{proj}_F(x)$ for all $t\geq 0$ and $x\in \mathbb{R}^n$ . In particular, $\varphi(t,x)$ varies along the segment from $x$ to $\mathrm{proj}_F(x)$.
\end{lemma}

\begin{proof}
Let $x \in \mathbb{R}$, let $x_\infty=\mathrm{proj}_F(x)$ and $\gamma_t=x_\infty+e^{-t}(x-x_\infty)$. It is immediate to check that $$\gamma_0=x$$ and that
\begin{equation}
\begin{split}
\frac{d}{dt} \gamma_t 
&= \frac{d}{dt} \big[x_\infty+e^{-t}(x-x_\infty)\big]\\ 
& = -e^{-t}(x-x_\infty) \\
& = -(\gamma_t-x_\infty).
\end{split}
\nonumber
\end{equation}
Moreover, let $t\geq 0$ and $z= \mathrm{proj}_F(\gamma_t)$. Therefore,
\begin{equation}
\begin{split}
\lVert x-x_\infty\rVert&=\lVert x-\gamma_t\rVert + \lVert \gamma_t - x_\infty\rVert\\
& \geq \lVert x-\gamma_t\rVert + \lVert \gamma_t - z\rVert\\
& \geq \lVert x-z\rVert.
\end{split}
\nonumber
\end{equation}

Hence $z=x_\infty$ by uniqueness of $\mathrm{proj}_F(x)$.
But then, $\gamma_t$ satisfies $\frac{d}{dt}\gamma_t=-\nabla d_F(\gamma_t)$ and therefore $\varphi(t,x)=\gamma_t$.
\end{proof}

Note that Lemma \ref{Cauchylip} implies that for every $x \in \mathbb{R}^n$, $\mathrm{proj}_F(x_\infty+s(x-x_\infty))=x_\infty$ for all $s \geq 0$, where $x_\infty = \mathrm{proj}_F(x)$. In fact, let $\psi(t,x)$ the flow of $\nabla d_F$, so that for all $t \geq 0$ and all $x \in \mathbb{R}^n$, $\varphi(t,\psi(t,x))=x$.
Applying Lemma \ref{Cauchylip}, we obtain $\psi(t,x)=x_\infty+e^t(x-x_\infty)$ and that $\mathrm{proj}_F(\psi(t,x))=x_\infty$. As a result, for all $s\geq 1$ $\mathrm{proj}_F(x_\infty+s(x-x_\infty))=x_\infty$. 

\noindent If $ 0 \leq s \leq 1$, then Lemma \ref{Cauchylip} already guarantees the previous equality.

We are now ready to prove Theorem \ref{C2dist}.
\\
\begin{proofof}{Theorem \ref{C2dist}}
    We first show that $F$ is convex.
    Let $y \in F$, $x \notin F$ and $x_\infty=\mathrm{proj}_F(x)$. By the above argument, for all $s\geq 0$
    $$\lVert x_\infty + s(x-x_\infty)-x_\infty\rVert^2 \leq \lVert x_\infty + s(x-x_\infty)-y\rVert^2,$$
 which can be rearranged as
    $$0 \leq 2s\langle x-x_\infty,x_\infty-y\rangle + \lVert x_\infty -y \rVert^2,$$
yielding $\langle x-x_\infty,y-x_\infty\rangle \leq 0$. Since $y$ is arbitrary, we have shown that every point $x \notin F$ can be separated from $F$ by an affine hyperplane. It follows that $F$ is convex.
Moreover, if $d_F$ is $C^2$ then $F$ is a convex subset of $\mathbb{R}^n$ that is a manifold without boundary by Lemma \ref{l3}. It is therefore an affine subspace of $\mathbb{R}^n$.
\end{proofof}

\bibliographystyle{alpha}
\bibliography{biblio.bib}

\end{document}